\documentclass[11pt, a4paper]{article}
\usepackage[cp1251]{inputenc}
\usepackage{amsmath} \usepackage{euscript}
\usepackage{amssymb}
\begin{document}

\newcounter{bnomer} \newcounter{snomer}
\renewcommand{\thesnomer}{\thebnomer.\arabic{snomer}}
\renewcommand{\refname}{\begin{center}\large{\textbf{References}}\end{center}}

\newcommand{\sect}[1]{%
\setcounter{snomer}{0} \refstepcounter{bnomer}
\begin{center}\large{\textbf{\arabic{bnomer}. {#1}}}\end{center}}
\newcommand{\defi}[1]{%
\refstepcounter{snomer}
\par\textbf{Definition \arabic{bnomer}.\arabic{snomer}. }{#1}\par}
\newcommand{\theo}[2]{%
\refstepcounter{snomer}
\par\textbf{Theorem \arabic{bnomer}.\arabic{snomer}. }{#2} {\emph{#1}} $\square$\par}
\newcommand{\mtheo}[1]{%
\refstepcounter{snomer}
\par\textbf{Theorem \arabic{bnomer}.\arabic{snomer}. }{\emph{#1}}\par}
\newcommand{\hypo}[1]{%
\refstepcounter{snomer}
\par\textbf{Conjecture \arabic{bnomer}.\arabic{snomer}. }{\emph{#1}}\par}
\newcommand{\theobp}[2]{%
\refstepcounter{snomer}
\par\textbf{Theorem \arabic{bnomer}.\arabic{snomer}. }{#2} {\emph{#1}}\par}
\newcommand{\theop}[2]{%
\refstepcounter{snomer}
\par\textbf{Theorem \arabic{bnomer}.\arabic{snomer}. }{\emph{#1}}
\par\textbf{Proof}. {#2} $\square$\par}
\newcommand{\exam}[1]{%
\refstepcounter{snomer}
\par\textbf{Example \arabic{bnomer}.\arabic{snomer}. }{#1}\par}
\newcommand{\deno}[1]{%
\refstepcounter{snomer}
\par\textbf{Definition \arabic{bnomer}.\arabic{snomer}. }{#1}\par}
\newcommand{\lemm}[1]{%
\refstepcounter{snomer}
\par\textbf{Lemma \arabic{bnomer}.\arabic{snomer}. }{\emph{#1}} $\square$\par}
\newcommand{\lemmp}[2]{%
\refstepcounter{snomer}
\par\textbf{Lemma \arabic{bnomer}.\arabic{snomer}. }{\emph{#1}}
\par\textbf{Proof}. {#2} $\square$\par}
\newcommand{\coro}[1]{%
\refstepcounter{snomer}
\par\textbf{Corollary \arabic{bnomer}.\arabic{snomer}. }{\emph{#1}} $\square$\par}
\newcommand{\mcoro}[1]{%
\refstepcounter{snomer}
\par\textbf{Corollary \arabic{bnomer}.\arabic{snomer}. }{\emph{#1}}\par}
\newcommand{\corop}[2]{%
\refstepcounter{snomer}
\par\textbf{Corollary \arabic{bnomer}.\arabic{snomer}. }{\emph{#1}}
\par\textbf{Proof}. {#2} $\square$\par}
\newcommand{\propp}[2]{%
\refstepcounter{snomer}
\par\textbf{Proposition \arabic{bnomer}.\arabic{snomer}. }{\emph{#1}}
\par\textbf{Proof}. {#2} $\square$\par}

\newcommand{\Ind}[3]{%
\mathrm{Ind}_{#1}^{#2}{#3}}
\newcommand{\Res}[3]{%
\mathrm{Res}_{#1}^{#2}{#3}}
\newcommand{\epsi}{\varepsilon}
\newcommand{\Supp}[1]{%
\mathrm{Supp}(#1)}

\newcommand{\reg}{\mathrm{reg}}
\newcommand{\sreg}{\mathrm{sreg}}
\newcommand{\codim}{\mathrm{codim}\,}
\newcommand{\chara}{\mathrm{char}\,}
\newcommand{\rk}{\mathrm{rk}\,}
\newcommand{\col}{\mathrm{col}}
\newcommand{\row}{\mathrm{row}}
\newcommand{\pho}{\hphantom{\quad}\vphantom{\mid}}
\newcommand{\wt}{\widetilde}
\newcommand{\ad}[1]{\mathrm{ad}_{#1}}

\newcommand{\vfi}{\varphi}
\newcommand{\lee}{\leqslant}
\newcommand{\gee}{\geqslant}
\newcommand{\Fp}{\mathbb{F}}
\newcommand{\Rp}{\mathbb{R}}
\newcommand{\Cp}{\mathbb{C}}
\newcommand{\ut}{\mathfrak{u}}
\newcommand{\at}{\mathfrak{a}}
\newcommand{\bt}{\mathfrak{b}}
\newcommand{\vt}{\mathfrak{v}}
\newcommand{\pt}{\mathfrak{p}}
\newcommand{\Po}{\EuScript{P}}
\newcommand{\Fo}{\EuScript{F}}
\newcommand{\Uo}{\EuScript{U}}
\newcommand{\Mo}{\mathcal{M}}
\newcommand{\Ro}{\mathcal{R}}
\newcommand{\Co}{\mathcal{C}}
\newcommand{\Lo}{\mathcal{L}}
\newcommand{\Ou}{\mathcal{O}}
\newcommand{\Sy}{\mathcal{Z}}
\newcommand{\Sb}{\mathcal{F}}
\newcommand{\Gr}{\mathcal{G}}

\author{Mikhail V. Ignatyev\thanks{Samara state university,
Department of algebra and geometry,\endgraf443011, ak. Pavlova, 1,
Samara, Russia, \texttt{mihail.ignatev@gmail.com}}}
\date{}
\title{Orthogonal subsets of classical root systems \\
and coadjoint orbits of unipotent groups} \maketitle

\sect{Introduction and statements of the main results}

Let $\Phi$ be a root system and $k$ an algebraic extension of a
finite field of sufficiently large {cha\-racte\-ris\-tic}~$p$. Let
$G$ be the classical matrix group over $k$ with the root system
$\Phi$, $U$ the subgroup of $G$ consists of all unipotent
lower-triangular matrices from $G$, $\Phi^+\subset\Phi$ the
corresponding set of~positive roots, and $\ut=\mathrm{Lie}(U)$ the
Lie algebra of $U$.

In the case $k=\Fp_q$ one can use the orbit method to describe
complex irreducible characters of~$U$ \cite{Kirillov1},
\cite{Kazhdan}: they are in one-to-one correspondence with the
orbits of the coadjoint representation of~$U$ in~the space $\ut^*$;
moreover, a lot of questions about representations can be
interpreted in terms of~orbits. Note that the problem of complete
description of orbits remains unsolved and seems to be very
difficult.

Let $D\subset\Phi^+$ be a subset consisting of pairwise orthogonal
roots. To each set of non-zero scalars $\xi=(\xi_{\beta})_{\beta\in
D}$ we assign the element of $\ut^*$ of the form
\begin{equation*}
f=f_{D, \xi}=\sum_{\beta\in D}\xi_{\beta}e_{\beta}^*
\end{equation*}
(by $e_{\beta}^*\in\ut^*$ we denote the covector dual to the root
vector $e_{\beta}\in\ut$ corresponding to a given root~$\beta$). By
$\Omega=\Omega_{D, \xi}$ we denote the orbit of $f$ under the
coadjoint action of $U$. We say that the orbit $\Omega$
is~\emph{associated} with the set $D$ and $f$ is the \emph{canonical
form} on this orbit.

The main goal of the paper is to compute the dimension of the orbit
$\Omega$ and to construct\linebreak a~polarization at $f$. (Recall
that a Lie subalgebra $\at\subset\ut$ is called a
\emph{polarization} of $\ut$ at a linear form $\lambda\in\ut^*$ if
$\lambda([\at,\at])=0$ and $\at$ is maximal among all subspaces of
$\ut$ with this property. Polarizations play an important role in
the explicit construction of the irreducible representation
corresponding\linebreak to~a~given orbit, see \cite[p. 274]{Kazhdan}
for the case $k=\Fp_q$.) As a consequence, we determine all possible
dimensions of irreducible representations of the group $U$.
Throughout the paper we suppose that $\Phi$ is of type $B_n$, $C_n$
or $D_n$ (the case of $A_n$ was considered by Alexander N. Panov in
\cite{Panov}). The paper generalizes results of~\cite{Ignatev2},
where those problems were solved by the author for the case
$\Phi=B_n,D_n$ and for orthogonal subsets of special kind.

\bigskip The paper is organized as follows. In section \ref{sect_pol}, we give nessesary
definitions. Then, for a given orthogonal subset $D$ we~construct
the subspace $\pt=\pt_D\subset\ut$  (see (\ref{formula_minus}) and
(\ref{formula_pol})). \mtheo{The subspace $\pt$ is a polarization of
$\ut$ at the form $f$.\label{mtheo_pol}}

In section \ref{sect_dim}, using the correspondence between the
dimensions of orbits and the codimensions of~polarizations \cite[p.
117]{Srinivasan} and induction by the rank of $\Phi$, we obtain a
formula for the dimension of~$\Omega$ (for the case of algebraically
closed field $k$). Precisely, let $W=W(\Phi)$ be the Weil group of
the root system $\Phi$ and $\sigma\in W$ the involution of the form
\begin{equation*}
\sigma=\prod_{\beta\in D}r_{\beta},
\end{equation*}
where $r_{\beta}$ is the reflection in the hyperplane orthogonal to
a given root $\beta$. Let $l(\sigma)$ be the length of $\sigma$ in
the Weil group, i.e, the length of the shortest (reduced)
representation of $\sigma$ as a product of simple reflections,
$s(\sigma)=|D|$, and $\vartheta$ the "defect"\ (see
(\ref{formula_defect})).

\mtheo{The dimension of $\Omega$ is equal to
$\dim\Omega=l(\sigma)-s(\sigma)-2\vartheta$.\label{mtheo_dim_orbit}}

In section \ref{sect_proofs}, using this theorem, we determine all
possible dimensions of irreducible representations of the group $U$
for the case of finite field $k$. Let $k=\Fp_q$ and $2\mu$ the
maximal possible dimension of~a~coadjoint orbit of $U$ (coadjoint
orbits are even dimensional). It was computed by
Carlos~A.~M.~Andr\`e and Ana M. Neto (see
\cite[Propositions~6.3,~6.6]{AndreNeto} and
(\ref{formula_max_dim})). \mcoro{The group $U$ has an irreducible
representation of dimension $N$ if and only if\linebreak$N=q^{l}$,
$0\lee l\lee\mu$.\label{mcoro_dim_rep}} (See \cite{Mukherjee} for
the case of $A_n$.) Section \ref{sect_proofs} also contains the
proofs of several technical results used in sections \ref{sect_pol},
\ref{sect_dim} and based on~detail (but elementary) studying of
roots from $D$.

\bigskip More generally, let $\Psi$ be an arbitrary root system and $D\subset\Psi^+$
an orthogonal subset. Let $\Uo$ be the maximal unipotent subgroup
of the Chevalley group over $k$ with the root system $\Psi$. As
above, let $\xi$ be the set of non-zero scalars from $k$ and
$\Omega=\Omega_{D,\xi}$ the coadjoint orbit of the group $\Uo$
associated with $D$.

\hypo{The dimension of the orbit $\Omega$ doesn't depend on $\xi$.
Further, $\dim\Omega\lee l(\sigma)-s(\sigma)$, where
$\sigma=\prod_{\beta\in D}r_{\beta}$ is the involution in the Weil
group $W(\Psi)$ corresponding to the subset $D$.} So in the paper we
prove this conjecture for classical root systems.

\bigskip The author is grateful to his scientific advisor professor Alexander
N. Panov for constant attention to this work.

\bigskip
\sect{A polarization at the form $f$}\label{sect_pol}

It's convenient to represent $\Phi$ as a subset of $\Rp^n$ (see
\cite{Bourbaki}): $\Phi=\pm\Phi^+$, where the set $\Phi^+$ of
positive roots has the form $\Phi^+=\Phi_0^+\cup\Phi_1^+$. Here
$\Phi_0^+=\{\epsi_i\pm\epsi_j,1\lee i<j\lee n\}$ and
\begin{equation*}
\Phi_1^+=\begin{cases} \varnothing,&\text{if }\Phi=D_n,\\
\{\epsi_i,1\lee i\lee n\},&\text{if }\Phi=B_n,\\
\{2\epsi_i,1\lee i\lee n\},&\text{if }\Phi=C_n
\end{cases}
\end{equation*}
($\{\epsi_i\}_{i=1}^n$ is the standard basis of $\Rp^n$).

Let $m=2n+1$ in the case $\Phi=B_n$ and $m=2n$ in the case
$\Phi=C_n$ or $D_n$. We'll index the rows and the columns of any
$m\times m$ matrix by the numbers $1,2,\ldots,n,0,-n,\ldots,-2,-1$
(if $m$ is even, then the index $0$ is omitted). We'll denote the
usual matrix units by $e_{a, b}$. By definition, $\ut$ is the
subalgebra of $\mathfrak{gl}_m(k)$ spanned by all $e_{\alpha}$,
$\alpha\in\Phi^+$, where
\begin{equation}
\begin{split}
&e_{\epsi_i-\epsi_j}=e_{j,i}-e_{-i,-j},\quad 1\lee i<j\lee n,\\
&e_{\epsi_i+\epsi_j}=e_{-j,i}-e_{-i,j},\quad 1\lee i<j\lee n,\\
&e_{\epsi_i}=e_{0, i}-e_{-i, 0},\quad e_{2\epsi_i}=e_{-i, i},\quad 1\lee i\lee n.\\
\end{split}\label{formula_ut}
\end{equation}

In the sequel, we assume that $\chara{k}$ is not less than $m$.
Under this assumption, the exponential map $\exp(x)=\sum_{i=0}^n
x^i/i!, x\in\ut$, is well-defined and bijective; moreover
$U=\exp(\ut)$ is a maximal unipotent subgroup of $G$ and
$\ut=\mathrm{Lie}(U)$. The group $U$ acts on $\ut$ via the adjoint
representation; the dual representation is called \emph{coadjoint}.
We'll denote the coadjoint action by $x.\lambda$, $x\in U$,
$\lambda\in\ut^*$.

\bigskip Now, let $D\subset\Phi$ be an \emph{orthogonal} subset
(i.e., subset consists of pairwise orthogonal roots),\linebreak
$\xi=(\xi_{\beta})_{\beta\in D}$ a set of non-zero scalars,
$\Omega=\Omega_{D,\xi}\subset\ut^*$ the associated coadjoint orbit,
and $f=f_{D,\xi}$ the canonical form on this orbit.

Let $i\neq j$. We assume without loss of generality that if
$\Phi=B_n$, then $D\cap\{\epsi_i,\epsi_j\}\leq1$, and if $\Phi=C_n$,
then $D\cap\{\epsi_i+\epsi_j,\epsi_i-\epsi_j\}\leq1$. Indeed, the
following proposition holds. \propp{$\mathrm{a)}$ Let
$\Phi=B_n$, $i<j$, $D\subset\Phi$ be an orthogonal subset containing
the roots $\epsi_i,\epsi_j$ and $\xi=(\xi_{\beta})_{\beta\in D}$ a
set of non-zero scalars. Let $D'=D\setminus\{\epsi_i\}$,
$\xi'=\xi\setminus\{\xi_{\epsi_j}\}$. Then
$\Omega_{D,\xi}=\Omega_{D',\xi'}$.\linebreak$\mathrm{b)}$ Let
$\Phi=C_n$, $D\subset\Phi$ be an orthogonal subset containing the
roots $\epsi_i-\epsi_j,\epsi_i+\epsi_j$,
$\xi=(\xi_{\beta})_{\beta\in D}$ be~a~set of non-zero scalars. Let
$D'=D\setminus\{\epsi_i-\epsi_j\}$,
$\xi'=\xi\setminus\{\xi_{\epsi_i-\epsi_j}\}$. Then
$\Omega_{D,\xi}=\Omega_{D',\xi'}$.\label{prop_min_D}}{a) Let
$f_{D',\xi'}$ be the canonical form on the orbit $\Omega_{D',\xi'}$
and $f'=\exp(ce_{\epsi_i-\epsi_j}).f_{D',\xi'}$ for some $c\in k^*$.
Then, by definition, for any $\alpha\in\Phi^+$
\begin{equation*}
f'(e_{\alpha})=f_{D',\xi'}(\exp\ad{-ce_{\epsi_i-\epsi_j}}(e_{\alpha}))=
f_{D',\xi'}(e_{\alpha})-c\cdot
f_{D',\xi'}(\ad{e_{\epsi_i-\epsi_j}}e_{\alpha})+
\dfrac{1}{2}c^2\cdot
f_{D',\xi'}(\ad{e_{\epsi_i-\epsi_j}}^2e_{\alpha})-\ldots
\end{equation*}

Recall that $[e_{\alpha},e_{\beta}]=c_{\alpha\beta}e_{\alpha+\beta}$
for all $\alpha,\beta\in\Phi^+$, and $c_{\alpha\beta}\neq0$ if and
only if $\alpha+\beta\in\Phi^+$. Since
$[e_{\epsi_i-\epsi_j},e_{\epsi_j}]=c_0e_{\epsi_i}$, $c_0\in k^*$, we
have $f'(e_{\epsi_j})= -\xi_{\epsi_i}\cdot c\cdot c_0$. On the other
hand, if $\alpha\neq\epsi_j$ and $f'(e_{\alpha})\neq0$, then there
exists $N\in\mathbb{Z}_{\gee0}$ such that
$\alpha+N(\epsi_i-\epsi_j)\in D\setminus\{\epsi_i,\epsi_j\}$. If
$N=0$, then $\alpha\in D'$. Suppose $N>0$. Then the inner products
$(\alpha,\epsi_i)$ and $(\alpha,\epsi_j)$ equal $-N$ and $N$
respectively, so $N = 1$ and~$\alpha=-\epsi_i+\epsi_j\notin\Phi^+$.
This stands in contradiction to the choice of $\alpha$.

Therefore, if $c=-\xi_{\epsi_j}/(c_0\cdot\xi_{\epsi_i})$, then
$f'=f$ and $\Omega_{D,\xi}=\Omega_{D',\xi'}$.

b) Let $f_{D',\xi'}$ be the canonical form on the orbit
$\Omega_{D',\xi'}$, and $f'=\exp(ce_{2\epsi_j}).f_{D',\xi'}$ for
some $c\in k^*$. Since
$[e_{2\epsi_j},e_{\epsi_i-\epsi_j}]=c_0e_{\epsi_i+\epsi_j}$, $c_0\in
k^*$, we have $f'(e_{\epsi_i-\epsi_j})= -\xi_{\epsi_i+\epsi_j}\cdot
c\cdot c_0$. On the other hand, if~$\alpha\neq\epsi_i-\epsi_j$ and
$f'(e_{\alpha})\neq0$, then there exists $N\in\mathbb{Z}_{\gee0}$
such that $\alpha+N\cdot2\epsi_j\in
D\setminus\{\epsi_i\pm\epsi_j\}$. If~$N=0$, then $\alpha\in D'$.
Suppose $N>0$. Then the inner products $(\alpha,\epsi_i-\epsi_j)$
and $(\alpha,\epsi_i+\epsi_j)$ equal $2N$ and $-2N$ respectively, so
$N=1$ and $\alpha=-2\epsi_j\notin\Phi^+$. This stands in
contradiction to the choice of $\alpha$.

Therefore, if
$c=-\xi_{\epsi_i-\epsi_j}/(c_0\cdot\xi_{\epsi_i+\epsi_j})$, then
$f'=f$ and $\Omega_{D,\xi}=\Omega_{D',\xi'}$. }

\bigskip The goal of this section is to construct a polarization of $\ut$ at $f$, i.e.,
to construct a subalgebra of $\ut$, which is a maximal $f$-isotropic
subspace. To do this, we need some more definitions.

According to (\ref{formula_ut}), we define the functions
\begin{equation*}
\begin{split}
&\col\colon\Phi^+\to\{1,\ldots,n\}\colon\col(\epsi_i\pm\epsi_j)=\col(\epsi_i)=\col(2\epsi_i)=i,\\
&\row\colon\Phi^+\to\{-n,\ldots,n\}\colon\row(\epsi_i\pm\epsi_j)=\mp j,\row(\epsi_i)=0,\row(2\epsi_i)=-i.\\
\end{split}
\end{equation*}
For an arbitrary $-n\lee i\lee n$ and $1\lee j\lee n$ the sets
\begin{equation*}
\begin{split}
&\Ro_i = \Ro_i(\Phi) = \{\alpha\in\Phi^+\mid \row(\alpha)=i\},\\
&\Co_j = \Co_j(\Phi) = \{\alpha\in\Phi^+\mid \col(\alpha)=j\}
\end{split}
\end{equation*}
are called the $i$\emph{th row} and the $j$\emph{th column} of
$\Phi^+$ respectively. Note that Proposition~\ref{prop_min_D}
implies ${|D\cap \Ro_i|\lee1}$ and $|D\cap \Co_j|\lee2$ for all $i,
j$ (furthermore, if $|D\cap \Co_j|=2$, then $D\cap
\Co_j=\{\epsi_j-\epsi_l,\epsi_j+\epsi_l\}$ and $\Phi=B_n$ or $D_n$).

\defi{Let $\beta\in\Phi^+$. Roots $\alpha, \gamma\in\Phi^+$ are called
$\beta$-\emph{singular} if their sum coincides with $\beta$. The set
of all $\beta$-singular roots is denoted by $S(\beta)$ (see
\cite{Andre3}, \cite{AndreNeto}, \cite{Mukherjee}).}

It's easy to see that singular roots have the following form:
\begin{equation}
\begin{split}
&S(\epsi_i-\epsi_j)=\bigcup_{l=i+1}^{j-1}\{\epsi_i-\epsi_l,\epsi_l-\epsi_j\},\quad
1\lee i < j\lee n,\\
&S(\epsi_i) = \bigcup_{l=i+1}^{n}\{\epsi_i-\epsi_l,\epsi_l\},\quad
S(2\epsi_i)=\bigcup_{l=i+1}^{n}\{\epsi_i-\epsi_l,\epsi_i+\epsi_l\},\quad
1\lee i\lee n,\\
&S(\epsi_i+\epsi_j)=\bigcup_{l=i+1}^{j-1}\{\epsi_i-\epsi_l,\epsi_l+\epsi_j\}
\cup\bigcup_{l=j+1}^n\{\epsi_i-\epsi_l,\epsi_j+\epsi_l\}\cup\\
&\bigcup_{l=j+1}^n\{\epsi_i+\epsi_l,\epsi_j-\epsi_l\}\cup S_{ij}
,\quad 1\lee i<j\lee n,\text{ where}
\end{split}\label{formula_sing_roots}
\end{equation}
\begin{equation*}
S_{ij}=\begin{cases}\{\epsi_i, \epsi_j\},&\text{if }\Phi=B_n,\\
\{\epsi_i-\epsi_j,2\epsi_j\},&\text{if }\Phi=C_n,\\
\varnothing,&\text{if }\Phi=D_n.
\end{cases}
\end{equation*}

For our purposes, it's convenient to partition the set $S(\beta)$
into two subsets $S^+(\beta)$ and $S^-(\beta)$, where
\begin{equation}
S^+(\beta)=\begin{cases} \{\epsi_i+\epsi_l,\quad i<l\lee
n\},&\text{if }\Phi=C_n\text{ and }\beta=2\epsi_i,\\
S(\beta)\cap \Co_{\col(\beta)}&\text{otherwise,}
\end{cases}\label{formula_sing_pm}
\end{equation}
and $S^-(\beta)=S(\beta)\setminus S^+(\beta)$ (note that
$S^+(\beta)\subset \Co_{\col(\beta)}$ for all $\beta$).

\deno{Let $j_1<\ldots<j_t$ be the numbers of columns containing
roots from $D$. Put\linebreak $\Mo=\Mo_D=\cup_{i=0}^t\Mo_{j_i}$,
where $j_0=0$, $\Mo_0=\varnothing$ and
\begin{equation}
\Mo_{j_i}=\{\gamma\in S^-(\beta)\mid\beta\in D\cap \Co_{j_i} \text{
and
}\gamma,\beta-\gamma\notin\cup_{l=0}^{i-1}\Mo_{j_l}\}\label{formula_minus}
\end{equation}}
\noindent for all $i=1,\ldots,t$. \exam{Let $\Phi=D_7$,
$D=\{\epsi_1\pm\epsi_5,\epsi_2\pm\epsi_6,\epsi_3+\epsi_4\}$. Then
$\Mo=\Mo_1\cup\Mo_2\cup\Mo_3$, where
$\Mo_1=\{\epsi_2\pm\epsi_5,\epsi_3\pm\epsi_5,\epsi_4\pm\epsi_5\}\cup
\Co_5$, $\Mo_2=\{\epsi_3\pm\epsi_6,\epsi_4\pm\epsi_6\}\cup \Co_6$
and $\Mo_3=\{\epsi_4\pm\epsi_7\}$.}

Now we'll define the subspace $\pt\subset\ut$ and prove that it's a
polarization of $\ut$ at the canonical form~$f$ on the orbit
$\Omega$. Namely, put $\Po=\Phi^+\setminus\Mo$ and
\begin{equation}
\pt=\pt_{D,\xi}=\sum_{\alpha\in\Po}ke_{\alpha}+\pt_0.\label{formula_pol}
\end{equation}
Here $\pt_0$ denotes the subspace constructed as follows. By
definition, it's spanned by all vectors $x$ of~the~form
$x=\xi_{\epsi_l+\epsi_j}\cdot e_{\epsi_l-\epsi_j}-
\xi_{\epsi_i-\epsi_j}\cdot e_{\epsi_l+\epsi_j}$, where
$\epsi_i-\epsi_j,\epsi_i+\epsi_j\in D$, $i<l<j$,
$\epsi_l-\epsi_j,\epsi_l+\epsi_j\in\Mo_i$ and~$D\cap
\Ro_{-l}=\varnothing$. In particular if $\Phi=C_n$, then $\pt_0=0$
for all $D$ (this follows from Proposition~\ref{prop_min_D}
b)).\exam{Let $\Phi$ and $D$ be as in the previous example. Then
$\pt_0$ is spanned by the vectors $\xi_{\epsi_1+\epsi_5}\cdot
e_{\epsi_2-\epsi_5}- \xi_{\epsi_1-\epsi_5}\cdot
e_{\epsi_2+\epsi_5}$, $\xi_{\epsi_1+\epsi_5}\cdot
e_{\epsi_3-\epsi_5}- \xi_{\epsi_1-\epsi_5}\cdot e_{\epsi_3+\epsi_5}$
and $\xi_{\epsi_2+\epsi_6}\cdot e_{\epsi_3-\epsi_6}-
\xi_{\epsi_2-\epsi_6}\cdot e_{\epsi_3+\epsi_6}$.}

\bigskip Now we'll prove two technical results which are also used in the next sections.
Let
\begin{equation}
\wt\Phi^+=\begin{cases} \Phi^+\setminus(\Co_1\cup \Ro_0),&\text{if
}\Phi=B_n\text{ and }D\cap \Co_1=\{\epsi_1\},\\
\Phi^+\setminus(\Co_1\cup \Co_j\cup \Ro_j\cup \Ro_{-j}),&\text{if
}D\cap \Co_1\neq\varnothing\text{ and }D\cap\Co_1\subset\{\epsi_1-\epsi_j,\epsi_i+\epsi_j\},\\
\Phi^+\setminus \Co_1&\text{otherwise}.\\
\end{cases}\label{forumula_wt_Phi}
\end{equation}
Notice that Proposition \ref{prop_min_D} a) implies $D=(D\cap
\Co_1)\cup(D\cap\wt\Phi^+)$.

Put $\wt\Phi=\pm\wt\Phi^+$. We note that, in fact, $\wt\Phi$ is
isomorphic to the root system of rank less then the~rank~of~$\Phi$.
Namely, let
\begin{equation}
\Phi'=\begin{cases}D_{n-1},&\text{if }\Phi=B_n,\text{ }D\cap \Co_1=\{\epsi_1\},\\
B_{n-2},C_{n-2},D_{n-2},&\text{if }\varnothing\neq D\cap
\Co_1\subset\Phi_0^+,\text{
where }\Phi=B_n,C_n,D_n\text{ resp.},\\
B_{n-1},C_{n-1},D_{n-1}&\text{otherwise, where
}\Phi=B_n,C_n,D_n\text{ resp.}
\end{cases}\label{formula_Phi_shtrih}
\end{equation}

\lemmp{There exists an isomorphism of root systems
$\wt\Phi\cong\Phi'$.}{It's enough to construct an one-to-one map
$\pi\colon\wt\Phi^+\to\Phi'^+$ that can be extended to an isometry
$\langle\wt\Phi^+\rangle_{\Rp}\to\langle\Phi'^+\rangle_{\Rp}$.

Let us define the number $m'$ for $\Phi'$ by the same rule as the
number $m$ for $\Phi$ (see the beginning of~the section). If $m'$ is
even, then put $n'= m'/2$, else put $n=(m'-1)/2$. Let's index the
columns of~roots from $\wt\Phi^+$ from $1$ to $n'$; let's index the
rows of these roots from $-n'$ to $n'$ (omitting the index $0$
in~the case of even $m'$). The required map $\pi$ is constructed.
\label{lemm_cong_Phi}}

We'll denote the isomorphism $\wt\ut\to\ut'$ that takes each
$e_{\alpha}$, $\alpha\in\wt\Phi^+$, to $e_{\pi(\alpha)}$ by the same
letter $\pi$. (Here
$\wt\ut=\sum_{\alpha\in\wt\Phi^+}ke_{\alpha}\subset\ut$ and
$\ut'=\sum_{\alpha\in\Phi'^+}ke_{\alpha}\subset\mathfrak{gl}_{m'}(k)$
is the Lie algebra of the maximal unipotent subgroup $U'$ of the
classical group $G'$ with the root system $\Phi'$.)

One can deduce from (\ref{formula_sing_roots}) and
(\ref{formula_sing_pm}) that if $\alpha+\gamma=\beta$ and $\alpha\in
S^+(\beta)$, then $\gamma\in S^-(\beta)$ (and vice versa). It's
straightforward to check that if $\gamma\in S^-(\beta)$, then
$\col(\gamma)\gee\col(\beta)$. Moreover, in this case
$\row(\gamma)=\row(\beta)$ and $\col(\gamma)=\row(\alpha)$, or
$\row(\gamma)=-\row(\alpha)$ and $\col(\gamma)=-\row(\beta)$ (here
$\alpha=\beta-\gamma$). \lemmp{$\mathrm{a)}$ Let $\beta\in D\cap
\Co_j$, $\alpha\in S^+(\beta)$ and $\beta-\alpha\in\Mo_j$. Then
$D\cap(\alpha+\Phi^+)\subset D\cap\Co_j$.\linebreak$\mathrm{b)}$ Let
$|D\cap\Co_j|=2$ (and so $\Phi\neq C_n$), $\gamma\in\Mo_j$ and
$\col(\gamma)\neq\pm\row(\beta)$, where
$D\cap\Co_j=\{\beta,\beta'\}$.\linebreak Then
$D\cap(\gamma+\Po)\subset D\cap\Co_j$.\label{lemm_one_sing}}{Let us
prove part a) (part b) can be proved similarly). Assume that there
exists $\beta'\in D\cap \Co_i$, $i\neq j$, such that
$\alpha+\delta=\beta'$ (the case $i=j$ is evident). Since $i\neq j$
and $S^+(\beta)\subset \Co_j$, $S^+(\beta')\subset \Co_i$
(see~(\ref{formula_sing_pm})), we conclude that $\alpha\in
S^-(\beta')$ and $\delta\in S^+(\beta')$.

Moreover, $j=\col(\beta)=\col(\alpha)\gee\col(\beta')=i$, so $j>i$.
But $\alpha\notin\Mo_i$ means that $\delta\in\Mo_s$ for some $s<i$
(see (\ref{formula_minus})). In particular there exists a root
$\beta''\in D\cap \Co_s$ such that $\delta\in S^-(\beta)$. Put
$\eta=\beta''-\delta$. If $\row(\delta)=-\row(\eta)$,
$\col(\delta)=-\row(\beta'')$, then
$i=\col(\beta')=\col(\delta)=-\row(\beta'')$, and the roots
$\beta',\beta''$ aren't orthogonal. This contradiction shows that
$\row(\delta)=\row(\beta'')$, $\col(\delta)=\row(\eta)$.

Similarly, if $\row(\alpha)=-\row(\delta)$,
$\col(\alpha)=\row(\beta')$, then
$j=\col(\beta)=\col(\alpha)=\row(\beta')$, and the~roots
$\beta,\beta'$ aren't orthogonal. This contradiction shows that
$\row(\alpha)=\row(\beta')$, $\col(\alpha)=\row(\delta)$. But~in
this case, $i=\col(\beta)=\col(\alpha)=\row(\delta)=\row(\beta'')$,
and the roots $\beta,\beta''$ aren't orthogonal. This~contradiction
proves the lemma. }

\bigskip Things are now ready to the proof of Theorem
\ref{mtheo_pol}. The proof immediately follows from the definition
of polarization and of two following Propositions.

\propp{The subspace $\pt$ is a subalgebra of
$\ut$.\label{prop_pol_subal}}{Denote $\ut_1=\sum_{\alpha\in
\Co_1}ke_{\alpha}$,
$\ut_2=\sum_{\alpha\in\Phi^+\setminus(\Co_1\cup\wt\Phi^+)}ke_{\alpha}$.
We see that $\pt=\pt_1+\pt_2+\wt\pt$ as vector spaces (here
$\pt_1=\pt\cap\ut_1$, $\pt_2=\pt\cap\ut_2$ and
$\wt\pt=\pt\cap\wt\ut$). The proof is by induction on the rank of
$\Phi$ (the~base can be checked directly). Let
$D'=\pi(D\cap\wt\Phi^+)$ and $\pt'$ be the subspace of $\ut'$
constructed by~the~rule (\ref{formula_pol}) applied to the
subset~$D'$ and the set of non-zero scalars $\xi'$ (this set
coincides with $\xi$ without scalars corresponding to the roots from
$D\cap \Co_1$; in particular if $D\cap\Co_1=\varnothing$, then
$\xi'=\xi$).

According to Lemma \ref{lemm_cong_Phi}, the rank of $\Phi'$ is less
than the rank of $\Phi$. Thus, by the inductive assumption, $\pt'$
is a subalgebra of $\ut'$. Hence, $\wt\pt$ is a subalgebra of
$\wt\ut$ (and of $\ut$) as the preimage of~a~subalgebra under the
morphism $\pi$. One can see that $\pt_1$ is a commutative ideal. So
it's enough to~prove that $[\pt_2+\wt\pt,\pt_2]\subset\pt$. By
definition, $\pt_2=\at+\bt$, where
$\at=\sum_{\alpha\in\Po\setminus(\Co_1\cup\wt\Phi^+)}ke_{\alpha}$,
and $\bt=\pt_2\cap\pt_0$. Consider the subspaces $\at$ and $\bt$ in
more details.

1. $[\pt_2+\wt\pt,\at]\subset\pt$. Indeed, if $\Phi=C_n$ and $D\cap
\Co_1=\{2\epsi_1\}$, or $\Phi=B_n$ and $D\cap \Co_1=\{\epsi_1\}$,
or~$D\cap \Co_1=\varnothing$, or $D\cap
\Co_1=\{\epsi_1-\epsi_j,\epsi_1+\epsi_j\}$ for some $j$, then
$\at=0$. Suppose that $D\cap \Co_1=\{\epsi_1-\epsi_j\}$.\linebreak
Then $\pt_1=\ut_1$, $\bt=0$ and $\at$ is spanned by the vectors
$e_{\alpha}$, $\alpha\in(\Ro_{-j}\cup
\Co_j)\setminus\{\epsi_1+\epsi_j\}$; in this case, the~inner
products $(\alpha, \epsi_j)$ are \emph{positive}. At the same time,
the roots $\epsi_i-\epsi_j$, $2\lee i\lee n$, belong to $\Mo_1$, so
if the coefficient of $e_{\gamma}$ in the sum $y=\sum
y_{\gamma}e_{\gamma}\in\pt_2+\wt\pt$ is non-zero, then the inner
products $(\gamma+\alpha,\epsi_j)$ are also positive. We conclude
that $\gamma+\alpha\in \Co_j\cup \Ro_{-j}$ and
$[e_{\gamma},e_{\alpha}]\in ke_{\gamma+\alpha}\subset\pt$. Since
$\gamma$ and $\alpha$ are arbitrary, $[y,\at]\subset\pt$ as
required.

Now, let $D\cap \Co_1=\{\epsi_1+\epsi_j\}$. Then $\pt_1=\ut_1$,
$\bt=0$ and $\at$ is spanned by the vectors
$e_{\alpha}$,\linebreak$\alpha\in
\Ro_j\setminus\{\epsi_1-\epsi_j\}$; in this case, the inner products
$(\alpha,\epsi_j)$ are \emph{negative}. At the same time, the roots
$\epsi_i+\epsi_j$, $2\lee i\lee n$, belong to $\Mo_1$ (as the root
$2\epsi_j$ in the case, $\Phi=C_n$), so if the coefficient of
$e_{\gamma}$ in the sum $y=\sum y_{\gamma}e_{\gamma}\in\pt_2+\wt\pt$
is non-zero, then the inner products $(\gamma+\alpha,\epsi_j)$ are
also negative. We~conclude that $\gamma+\alpha\in\Ro_j$ and
$[e_{\gamma},e_{\alpha}]\in ke_{\gamma+\alpha}\subset\pt$. Since
$\gamma$ and $\alpha$ are arbitrary, $[y,\at]\subset\pt$ as
required.

2. $[\pt_2+\wt\pt,\bt]\subset\pt$. This follows from Lemma
\ref{lemm_p_0_1}. }

\bigskip
\propp{The subspace $\pt$ is a maximal $f$-isotropic
subspace.\label{prop_polar}}{1. First, let us prove that $\pt$ is an
$f$-isotropic subspace. Suppose that
$y,z\in\sum_{\alpha\in\Po}ke_{\alpha}$ (recall that
$\Po=\Phi^+\setminus\Mo$). In this case, $[y,
z]\in\sum_{\alpha,\gamma\in\Po}ke_{\alpha+\gamma}$ (we assume
$e_{\alpha}=0$, if $\alpha\notin\Phi^+$). It follows from $f([y,
z])\neq0$ that there exist $\alpha,\gamma\in\Po$ such that
$\alpha+\gamma\in D$. But this stands in contradiction with the
definition of $\Mo$ (see (\ref{formula_minus})). Indeed, the set
$\Mo$ contains either \emph{one} or \emph{two} of roots from each
pair of~$\beta$-singular roots which sum equals $\beta\in D$. Thus,
$\Po$ cannot contain the roots $\alpha,\gamma$ at the same time.

Now, let $x=\xi_{\epsi_i+\epsi_j}\cdot e_{\epsi_l-\epsi_j}-
\xi_{\epsi_i-\epsi_j}\cdot e_{\epsi_l+\epsi_j}\in\pt_0$, $i<l<j$ (so
$\Phi\neq C_n$ and $\beta,\beta'\in D$, where
$\beta=\epsi_i-\epsi_j$, $\beta'=\epsi_i+\epsi_j$). If
$\alpha\in\Po$, then, according to Lemma \ref{lemm_one_sing} b),
$\alpha+(\epsi_l\pm\epsi_j)\in D$ implies $\alpha=\epsi_i-\epsi_l$.
But in this case, $f([x, e_{\alpha}])=0$. On the other hand, if
$x'\in\pt_0$ and $[x,x']\neq0$, then, obviously,
$x'=\xi_{\epsi_i+\epsi_j}\cdot e_{\epsi_s-\epsi_j}-
\xi_{\epsi_i-\epsi_j}\cdot e_{\epsi_s+\epsi_j}$, $i<s<j$. Assume
that, for instance, $s<l$. Then $[x, x']\in ke_{\epsi_s+\epsi_l}$.
But if $\epsi_s+\epsi_l\in D$, then $x'$ \emph{cannot} belongs to
$\pt_0$ by definition of this space (see (\ref{formula_pol})).

2. Let us now show that $\pt$ is maximal (with respect to the
inclusion order) among all $f$-isotropic subspaces. Suppose that
there exists $y\notin\pt$ such that $\pt+ky$ is an isotropic
subspace. Let
\begin{equation*}
y=\sum_{\gamma\in\Mo}y_{\gamma}e_{\gamma},\quad y_{\gamma}\in k^*.
\end{equation*}

Pick a root $\gamma_0$ such that $y_{\gamma_0}\neq0$; by definition,
$\gamma_0\in\Mo_i$ for some $i$. In other words, there exist
$\beta\in D\cap \Co_i$ and $\alpha_0\in S^+(\beta)$ such that
$\beta=\alpha_0+\gamma_0$, and $\alpha_0\in\Po$, i.e.,
$e_{\alpha_0}\in\pt$. Hence, $[e_{\alpha_0},e_{\gamma_0}]=c\cdot
e_{\beta}$, $c\in k^*$. Applying Lemma \ref{lemm_one_sing} a), we
see that if $D\cap \Co_i=\{\beta\}$, then $\alpha_0+\gamma\notin D$
for all $\gamma\neq\gamma_0$, Thus,
\begin{equation*}
f([e_{\alpha_0},y])=f([e_{\alpha_0},y_{\gamma_0}e_{\gamma_0}])=y_{\gamma_0}\cdot
c\cdot\xi_{\beta}\neq0.
\end{equation*}

Lemma \ref{lemm_p_0_2} guarantees that if $|D\cap \Co_i|=2$, then
there exists $x'\in\pt+ky$ such that $f[y,x']\neq0$. We~see that
$\pt$ \emph{can't} be included into an isotropic subspace of higher
dimension. The result follows. }

\bigskip The proof of Theorem \ref{mtheo_pol} is complete. In some cases
(for example, if $\Phi=C_n$) one can use it to~compute the dimension
of an orbit associated with an orthogonal subset. (From now on to
the end of~the next section, we assume that the ground field $k$ is
algebraically closed.) \corop{Suppose $|D\cap \Co_j|\lee1$ for all
$1\lee j\lee n$. Then $\dim\Omega=2\cdot|\Mo|$.}{Indeed, the
dimension of an orbit is twice to the codimension of a polarization
at an arbitrary point on this orbit \cite[p. 117]{Srinivasan}. But
in our case, the codimension of $\pt$ equals $|\Mo|$, because
$\pt_0=0$. }

However, we'll obtain an explicit formula for the dimension of an
orbit associated with an \emph{arbitrary} orthogonal subset (see
Theorem \ref{mtheo_dim_orbit}). In order to prove this formula we'll
consider the involution in~the~Weil group that equals to the product
of reflections corresponding to the roots from $D$.

\bigskip
\sect{The dimension of the orbit $\Omega$}\label{sect_dim}

Let $D$, $\xi$, $\Omega$, $\pt$ be as in the previous section.
Recall that we defined the root system $\Phi'$ of rank less than the
rank of $\Phi$ (see (\ref{formula_Phi_shtrih})) and the subset
$\wt\Phi^+\subset\Phi^+$. We also constructed the (one-to-one) map
$\pi\colon\wt\Phi^+\to\Phi'^+$, which can be extended to the
isomorphism of root systems (see Lemma \ref{lemm_cong_Phi}), and put
$D'=\pi(D\cap\wt\Phi^+)$. Finally, we defined the subalgebra
$\pt'\subset\ut'$ (see the proof of Proposition
\ref{prop_pol_subal}). Notice that $D=(D\cap \Co_1)\cup\pi^{-1}(D')$
and these subsets are disjoint.

For simplicity, denote
\begin{equation}
r=\begin{cases} |\Co_1|+|S^-(\epsi_1\mp\epsi_j)|,&\text{if }D\cap
\Co_1=\{\epsi_1\pm\epsi_j\},\\
|\Co_1|+\#\{l\mid1<l<j\text{ and }D\cap
\Ro_{-l}=\varnothing\},&\text{if
}D\cap \Co_1=\{\epsi_1-\epsi_j,\epsi_1+\epsi_j\},\\
|\Co_1\cap\Po|&\text{otherwise}.
\end{cases}\label{formula_r}
\end{equation}
\lemmp{The dimensions of $\pt$ and $\pt'$ satisfy the equality
$\dim\pt=\dim\pt'+r$.\label{lemm_pt_pt_shtrih}}{One can represent
$\pt$ as a sum of vector spaces $\pt=\pt_1+\pt_2+\wt\pt$, where
$\pt_i=\pt\cap\ut_i$, $i=1,2$, $\ut_1=\sum_{\alpha\in
\Co_1\cap\Po}ke_{\alpha}$,
$\ut_2=\sum_{\alpha\in\Phi^+\setminus(\Co_1\cup\wt\Phi^+)}ke_{\alpha}$
and $\wt\pt=\pt\cap\sum_{\alpha\in\wt\Phi^+}ke_{\alpha}$ (see the
proof of Proposition~\ref{prop_pol_subal}). Since $\wt\pt\cong\pt'$,
we obtain $\dim\pt-\dim\pt'=\dim\pt_1+\dim\pt_2=|\Co_1\cap\Po|+
|\Po\setminus(\Co_1\cup\wt\Phi^+)|+\dim(\pt_0\cap\ut_2)$. It's
straightforward to check that the RHS of the last formula equals
$r$. }

Let $W$ be the Weil group of the root system $\Phi$. For an
arbitrary $\alpha\in\Phi^+$, by $r_{\alpha}\in W$ we denote the
reflection on the hyperplane orthogonal to $\alpha$. Consider the
following involution (i.e., the element of order two) in $W$:
\begin{equation*}
\sigma=\sigma_D=\prod_{\beta\in D}r_{\beta}
\end{equation*}
(commuting reflections $r_{\beta}$ are taken in any fixed order). We
define the involution $\sigma'$ in the Weil group~$W'$ of the root
system $\Phi'$ similarly (starting from the subset
$D'\subset\Phi'^+$).

By $l(\sigma)$ we denote the length of the shortest (reduced)
representation of $\sigma$ as a product of simple reflections. (In
other words, $l(\sigma)$ is the length of $\sigma$ as an element of
the Weil group). Let $s(\sigma)=|D|$. Let $l'(\sigma')$ and
$s'(\sigma')$ be defined by the similar rule.

In order to describe the dimension of an orbit $\Omega$ associated
with the orthogonal subset $D$, we'll define the number $\vartheta$.
By definition,
\begin{equation}
\begin{split}
&\vartheta=d_1+d_2+d_3+d_4,\text{ where}\\
&d_1=\#\{(i, j, l, s)\mid
i<l<s<j\text{ and }
\epsi_i-\epsi_j,\epsi_i+\epsi_j,\epsi_l+\epsi_s\in D\},\\
&d_2=\#\{(i, j, l, s)\mid i<l<j<s\text{ and }
\epsi_i-\epsi_j,\epsi_i+\epsi_j,\epsi_l-\epsi_s,\epsi_l+\epsi_s\in D\},\\
&d_3=\#\{(i, j)\mid\epsi_i+\epsi_j\in D\text{ and }i>l,\text{ where
}
D\cap \Ro_0=\{\epsi_l\}\},\\
&d_4=\#\{(i, j)\mid\epsi_i-\epsi_j,\epsi_i+\epsi_j\in D\text{ and }
i<j<l,\text{ where }D\cap \Ro_0=\{\epsi_l\}\}.\\
\end{split}\label{formula_defect}
\end{equation}
Note that if $\Phi=C_n$, then $\vartheta=0$ for all
$D\subset\Phi^+$. If $d_3\neq0$ or $d_4\neq0$, then $\Phi=B_n$ and
$D\cap \Ro_0\neq\varnothing$; in this case, $|D\cap \Ro_0|=1$ (see
Proposition \ref{prop_min_D} a)), so $d_3$ and $d_4$ are
well-defined.

We define the number $\vartheta'$ similarly (starting from the
subset $D'\subset\Phi'^+$). Obviously,
\begin{equation*}
s(\sigma)=s'(\sigma')+|D\cap\Co_1|,
\end{equation*}so the proof of
Theorem \ref{mtheo_dim_orbit} is based on comparing $l(\sigma)$ with
$l'(\sigma')$ and $\vartheta$ with $\vartheta'$ resp.

For a given involution $\tau\in W$, by $\Phi_{\tau}$ we denote the
set of positive roots such that their images under the action of
$\tau$ are negative:
$\Phi_{\tau}=\{\alpha\in\Phi^+\mid\tau(\alpha)\in\Phi^-\}$. It's
well-known that $l(\tau)=|\Phi_{\tau}|$, so~we are to compare the
numbers of elements of the sets $\Phi_{\sigma}$ and
$\Phi'_{\sigma'}$. Let $\wt D=D\cap\wt\Phi^+=\pi^{-1}(D')$, and
$\wt\sigma\in W$ be the involution of the form $\prod_{\beta\in\wt
D}r_{\beta}$.

The intersection of $\wt\Phi^+$ with the first column of $\Phi^+$ is
empty. Similarly, if $\alpha\in D\cap\Co_1$, then
$\wt\Phi^+\cap\row(\alpha)=\wt\Phi^+\cap\col(\alpha)=\varnothing$.
Thus, $\sigma(\alpha)=\wt\sigma(\alpha)$ for all
$\alpha\in\wt\Phi^+$. Therefore,
$\Phi_{\sigma}\cap\wt\Phi^+=\pi^{-1}(\Phi'_{\sigma'})$ and
$|\Phi_{\sigma}\cap\wt\Phi^+|=|\Phi'_{\sigma'}|=l'(\sigma')$. So it
remains to study the action of $\sigma$ on
$\Phi^+\setminus\wt\Phi^+$.

\bigskip\lemmp{Suppose $D\cap \Co_1=\{\epsi_1-\epsi_j\}$. Then
$l(\sigma)=l'(\sigma')+|S(\epsi_1-\epsi_j)|+1$.\label{lemm_i_m_j}}{In
our case, $\Phi^+\setminus\wt\Phi^+=\Co_1\cup \Co_j\cup \Ro_j\cup
\Ro_{-j}$ (see (\ref{forumula_wt_Phi})). Here $\Phi_{\sigma}\cap
\Co_1=S^+(\epsi_1-\epsi_j)\cup\{\epsi_1-\epsi_j\}$.\linebreak
If~${\alpha\in \Co_j\cup \Ro_{-j}}$, then the inner product
$(\alpha,\epsi_j)$ is positive. The Weil group acts by orthogonal
{trans\-for\-ma\-tions}, so $(\sigma(\alpha),\epsi_1)>0$. It follows
that $\sigma(\alpha)>0$ (i.e., belongs to $\Phi^+$). If $\alpha\in
\Ro_j\setminus\{\epsi_1-\epsi_j\}=S^-(\epsi_1-\epsi_j)$, then
$(\alpha,\epsi_j)<0$. Thus, $(\sigma(\alpha),\epsi_1)<0$ and so
$\sigma(\alpha)<0$ (i.e., belongs to $\Phi^-$). Hence,
\begin{equation*}
l(\sigma)=|\Phi'_{\sigma'}|+|S^+(\epsi_1-\epsi_j)| +
1+|S^-(\epsi_1-\epsi_j)|= l'(\sigma')+|S(\epsi_1-\epsi_j)|+1
\end{equation*}
as required. }

\lemmp{Suppose $D\cap \Co_1=\{\epsi_1+\epsi_j\}$. Then
$l(\sigma)=l'(\sigma')+|S(\epsi_1+\epsi_j)|+1$.\label{lemm_i_p_j}}{As
in the previous Lemma, $\Phi^+\setminus\wt\Phi^+=\Co_1\cup \Co_j\cup
\Ro_j\cup \Ro_{-j}$ (see (\ref{forumula_wt_Phi})). Here
$\Phi_{\sigma}\cap
\Co_1=\{\epsi_1\pm\epsi_i,i<j\}\cup\{\epsi_1+\epsi_j\}\cup S_1$ (if
$\Phi=D_n$, then $S_1$ is empty; if $\Phi=B_n$, then
$S_1=\{\epsi_1\}$; if $\Phi=C_n$, then $S_1=\{2\epsi_1$\}). By the
way, $\Phi_{\sigma}\cap \Co_1$ consists of
$|S^+(\epsi_1+\epsi_j)|+1$ roots.

If $\alpha\in
(\Co_j\cup\Ro_{-j})\setminus\{\epsi_1+\epsi_j\}=S^-(\epsi_1+\epsi_j)$,
then $(\alpha,\epsi_j)>0$, so $(\sigma(\alpha),\epsi_1)<0$ and
$\sigma(\alpha)<0$. If~$\alpha\in \Ro_j$, then $(\alpha,\epsi_j)<0$,
$(\sigma(\alpha),\epsi_1)>0$ and $\sigma(\alpha)>0$. Hence,
\begin{equation*}
\begin{split}l(\sigma)&=|\Phi'_{\sigma'}|+|S^+(\epsi_1+\epsi_j)|+1+
|S^-(\epsi_1+\epsi_j)|=l'(\sigma')+|S(\epsi_1+\epsi_j)|+1
\end{split}
\end{equation*}
as required. }

\lemmp{$\mathrm{a)}$ Suppose $\Phi=B_n$ and $D\cap
\Co_1=\{\epsi_1\}$. Then
$l(\sigma)=l'(\sigma')+|\Co_1|+2\cdot\#\{\beta\in\wt
D\mid\row(\beta)<0\}$. $\mathrm{b)}$ Suppose $\Phi=C_n$ and $D\cap
\Co_1=2\epsi_1$. Then
$l(\sigma)=l'(\sigma')+|\Co_1|$.\label{lemm_i_2i}}{a) In this case,
$\Phi^+\setminus\wt\Phi^+=\Co_1\cup \Ro_0$. It's clear that
$\Phi_{\sigma}\cap \Co_1=\Co_1$. If $\alpha\in\Ro_0$, then
$\sigma(\alpha)=\wt\sigma(\alpha)=\pm\epsi_l$ for some $l$. Since
$\epsi_1$ is orthogonal to all other roots from $\Ro_0$, we see that
$\sigma(\alpha)<0$ if and only if $\epsi_i+\epsi_l\in\wt D$. Hence,
$l(\sigma)=|\Phi'_{\sigma'}|+|\Co_1|+2\cdot\#\{\beta\in\wt
D\mid\row(\beta)<0\}$ as required.

b) Evident: $\Phi^+\setminus\wt\Phi^+=\Co_1$ is contained in
$\Phi_{\sigma}$. }

The case $|D\cap \Co_1|=2$ is considered in Lemma \ref{lemm_i_pm_j}.

\bigskip The following Proposition plays the key role in the proof of Theorem
\ref{mtheo_dim_orbit}. \propp{Let $D,D',\sigma,\sigma'$ and $r$ be
as above. Then
\begin{equation*}
l(\sigma)-s(\sigma)-2\vartheta=l'(\sigma')-s'(\sigma')-2\vartheta'+
2(|\Phi^+\setminus\wt\Phi^+|-r).
\end{equation*}
\label{prop_Fo_Fo_shtrih}}{For simplicity, denote
$\Fo=l(\sigma)-s(\sigma)-2\vartheta$ and
$\Fo'=l'(\sigma')-s'(\sigma')-2\vartheta'$. The proof is by
consideration of different variants of $D\cap \Co_1$.

1. $D\cap \Co_1=\varnothing$. Of course, in this case,
$l(\sigma)=l'(\sigma')$, $s(\sigma)=s'(\sigma')$ and
$\vartheta=\vartheta'$. On the other hand,
$\Phi^+\setminus\wt\Phi^+=\Co_1$ (see (\ref{forumula_wt_Phi})) and
$r=|\Co_1|$ (see (\ref{formula_r})) as required.

2. $D\cap \Co_1=\{\epsi_1\pm\epsi_j\}$. Here
$s(\sigma)=s'(\sigma')$, $\vartheta=\vartheta'$ and, according to
Lemmas \ref{lemm_i_m_j}, \ref{lemm_i_p_j},
$l(\sigma)=l'(\sigma')+|S(\epsi_1\pm\epsi_j)|+1$. So,
$\Fo-\Fo'=|S(\epsi_1\pm\epsi_j)|$. At the same time, using
(\ref{forumula_wt_Phi}), (\ref{formula_r}) and the fact that
$\Co_1\cup \Co_j\cup \Ro_j\cup \Ro_{-j}=\Co_1\cup
S^-(\epsi_1-\epsi_j)\cup S^-(\epsi_1+\epsi_j)$, we obtain
\begin{equation*}
|\Phi^+\setminus\wt\Phi^+|-r=|\Co_1\cup \Co_j\cup \Ro_j\cup
\Ro_{-j}|-(|S^-(\epsi_1\mp\epsi_j)|+|\Co_1|)=|S^-(\epsi_1\pm\epsi_j)|.
\end{equation*}
The last number is two times less than $|S(\epsi_1\pm\epsi_j)|$ as
required.

3. $D\cap \Co_1=\{\epsi_1\}$ ($\Phi=B_n$). By
(\ref{forumula_wt_Phi}), (\ref{formula_r}), (\ref{formula_defect})
and Lemma \ref{lemm_i_2i} a), we have
$|\Phi^+\setminus\wt\Phi^+|-r=(|\Co_1|+|S^-(\epsi_1)|)-|\Co_1|=|S^-(\epsi_1)|=n-1$.
At the same time,
$l(\sigma)=l'(\sigma')+|\Co_1|+2\cdot\#\{\beta\in\wt
D\mid\row(\beta)<0\}$, $s(\sigma)=s'(\sigma')+1$
and~$\vartheta=\vartheta'+\#\{(i, j)\mid\epsi_i+\epsi_j\in D\text{ and
}i>1\}=\vartheta'+\#\{\beta\in\wt D\mid\row(\beta)<0\}$,
so~$\Fo-\Fo'=|\Co_1|-1=(2n-1)-1=2(n-1)$ as required.

4. $D\cap \Co_1=\{2\epsi_1\}$ ($\Phi=C_n$). By
(\ref{forumula_wt_Phi}), (\ref{formula_r}), (\ref{formula_defect})
and Lemma \ref{lemm_i_2i} b), we have
$|\Phi^+\setminus\wt\Phi^+|-r=|\Co_1|-|\Co_1\cap\Po|=(2n-1)-n=n-1$.
At the same time, $l(\sigma)=l'(\sigma')+|\Co_1|$,
$s(\sigma)=s'(\sigma')+1$ and~$\vartheta=\vartheta'=0$, so
$\Fo-\Fo'=|\Co_1|-1=(2n-1)-1=2(n-1)$ as required.

5. $D\cap \Co_1=\{\epsi_1-\epsi_j,\epsi_1+\epsi_j\}$ ($\Phi=B_n$ or
$D_n$). By (\ref{forumula_wt_Phi}) and (\ref{formula_r}), we get
\begin{equation*}
\begin{split}|\Phi^+\setminus\wt\Phi^+|-r&=|\Co_1\cup \Co_j\cup \Ro_j\cup
\Ro_{-j}|-|\Co_1|-\#\{l\mid1<l<j\text{ and }D\cap
\Ro_{-l}=\varnothing\}=\\
&=m-4-\#\{l\mid1<l<j\text{ and }D\cap \Ro_{-l}=\varnothing\}=\\
&=m-4-\#\{l\mid1<l<j\text{ and }\wt D\cap \Ro_{-l}=\varnothing\}=\\
&=m-4-(j-2)+\#\{l\mid1<l<j\text{ and }\wt D\cap
\Ro_{-l}\neq\varnothing\}=\\
&=m-j-2+\#\{l\mid1<l<j\text{ and }\wt D\cap
\Ro_{-l}\neq\varnothing\}.
\end{split}
\end{equation*}
On the other hand, $s(\sigma)=s'(\sigma')+2$. Comparing
(\ref{formula_defect}) with (\ref{formula_i_pm_j}) (see Lemma
\ref{lemm_i_pm_j}), we obtain
\begin{equation*}
\begin{split}
\Fo-\Fo'&=|\Co_1|+|\Co_j|-2+2\cdot\#\{(l,s)\mid1<l<s<j\text{ and
}\epsi_l+\epsi_s\in\wt D\}=\\
&=(m-2)+(m-2j)-2+2\cdot\#\{(l,s)\mid1<l<s<j\text{ and
}\epsi_l+\epsi_s\in\wt D\}=\\
&=2(m-j-2)+2\cdot\#\{s\mid1<s<j\text{ and }\wt D\cap
\Ro_{-s}\neq\varnothing\}.
\end{split}
\end{equation*}
To conclude the proof, it remains to replace $s$ by $l$ in the last
formula. }

\bigskip Combining this Proposition with Lemma
\ref{lemm_pt_pt_shtrih}, we'll now prove Theorem
\ref{mtheo_dim_orbit}.

\textbf{Proof of Theorem \ref{mtheo_dim_orbit}.} The proof is by
induction on the rank of $\Phi$ (the base is checked directly). Let
$\Omega'=\Omega_{D',\xi'}\subset\ut'^*$ be the orbit of the element
$f_{D',\xi'}$ under the coadjoint action of the group $U'$. By~the
inductive assumption,
$\dim\Omega'=l(\sigma')-s(\sigma')-2\vartheta'$. Since the dimension
of an orbit is twice to~the codimension of a polarization at a point
on this orbit \cite[p. 117]{Srinivasan}, we deduce from
Proposition~\ref{prop_Fo_Fo_shtrih} and Lemma
\ref{lemm_pt_pt_shtrih} that the dimension of the orbit $\Omega$
equals
\begin{equation*}
\begin{split}
&2\cdot\codim\pt=2(|\Phi^+|-\dim\pt)=2(|\Phi^+|-\dim\pt'-r)=\\
&=2(|\Phi^+|-r-|\Phi'^+|+\codim\pt')=2(|\Phi^+\setminus\wt\Phi^+|-r)+\dim\Omega'=\\
&=l(\sigma)-s(\sigma)-2\vartheta-(l'(\sigma)-s'(\sigma')-2\vartheta')+l'(\sigma)-s'(\sigma')-2\vartheta'=\\
&=l(\sigma)-s(\sigma)-2\vartheta\quad\text{as required. }\square
\end{split}
\end{equation*}
\exam{Let $\Phi=B_7$ and
$D=\{\epsi_1-\epsi_6,\epsi_1+\epsi_6,\epsi_2,\epsi_3-\epsi_7,\epsi_3+\epsi_7,\epsi_4+\epsi_5\}$.
We see that $s(\sigma)=|D|=6$, $l(\sigma)=|\Phi_{\sigma}|=48$ (one
can find $\Phi_{\sigma}$ explicitly). On the other hand, by
(\ref{formula_defect}) we~obtain $d_1=\#\{(1, 6, 4, 5), (3, 7, 4,
5)\}=2$, $d_2=\#\{(1, 6, 3, 7)\}=1$, $d_3=\#\{(3, 7), (4, 5)\}=2$
and $d_4=\#\{(1, 6)\}=1$ (because $D\cap\Ro_0=\{\epsi_2\}$). Hence,
$\vartheta=d_1+d_2+d_3+d_4=6$ and
\begin{equation*}
\dim\Omega=l(\sigma)-s(\sigma)-2\vartheta=48-6-12=30.
\end{equation*}}

\bigskip \sect{Dimensions of representations of $U$ and proofs}\label{sect_proofs}

From now on, let $k=\Fp_q$ be a finite field with $q$ elements (so
$U$ be a finite group). Using Theorem~\ref{mtheo_dim_orbit} and the
correspondence between irreducible finite-dimensional complex
representations of $U$ and coadjoint orbits, we'll now describe all
possible dimensions of these representations. Let $K$
be~the~algebraic closure of the field $k$, $\ut_K$ the subalgebra of
$\mathfrak{gl}_m(K)$ spanned by vectors
of~the~form~(\ref{formula_ut}), and $U_K=\exp(\ut_K)$. If
$f\in\ut^*\subset\ut_K^*$, then by $\Omega\subset\ut^*$ (resp.
$\Omega_K\subset\ut_K^*$) we'll denote its orbit under the~coadjoint
action of the group $U$ (resp. of the group $U_K$).

According to \cite[Proposition 2]{Kazhdan}, there is a one-to-one
correspondence between coadjoint orbits of $U$ and classes of
isomorphic irreducible representations of $U$; moreover, if the
orbit $\Omega$ corresponds to~a~given representation $V$, then $\dim
V=\sqrt{\Omega}=q^{\dim\Omega_K/2}$. Let
\begin{equation}
\mu=\begin{cases}n(n-1)/2,&\text{if }\Phi=B_n\text{ or }C_n,\\
n(n-1)/2,&\text{if }\Phi=D_n\text{ and }n\text{ is even},\\
(n-1)^2/2,&\text{if }\Phi=D_n\text{ and }n\text{ is odd}.
\end{cases}\label{formula_max_dim}
\end{equation}

If an orbit $\Omega_K$ is of maximal dimension, then its dimension
equals $2\mu$ \cite[Pro\-po\-si\-tions~6.3,~6.6]{AndreNeto}.
Corollary \ref{mcoro_dim_rep} claims that there exists a
representation of the group $U$ of dimension $N$ if and only if
$N=q^l$, where $0\lee l\lee\mu$. To prove this, it remains to find
an orbit $\Omega_K$ of dimension $2l$. Theorem \ref{mtheo_dim_orbit}
shows that it's enough to construct an orthogonal subset
$D\subset\Phi^+$ such that $l(\sigma)-s(\sigma)-2\vartheta=2l$
(in~fact, we'll deal with subsets such that $\vartheta=0$).

\textbf{Proof of Corollary \ref{mcoro_dim_rep}.} For an arbitrary
$1\lee j\lee [n/2]$, set $\beta_j=\epsi_{2j-1}+\epsi_{2j}$ and
$s_j=|S^+(\beta_j)|$. It's easy to check that $s_1+\ldots+s_t=\mu$,
where $t=[n/2]$ for $\Phi=B_n$ or $C_n$, and $t=[(n-1)/2]$ for
$\Phi=D_n$ (see (\ref{formula_sing_roots}) and
(\ref{formula_sing_pm})). We note also that if $\alpha\in\Co_{2j-1}$
and $\row(\alpha)$ runs $2j, 2j+1\ldots,n,0,-n,\ldots,-2j+1,-2j$
(the~index $0$ is omitted for even $m$), then $|S^+(\alpha)|$ runs
$0,1,\ldots,s_j$ respectively.

Let $0\lee l\lee\mu$. If $l\lee s_1$, then, as mentioned above,
there exists $\beta\in \Co_1$ such that $|S^+(\beta)|=l$.
Let~$D=\{\beta\}$, then $\Phi_{\sigma}=\Phi_{r_{\beta}}$. But for an
arbitrary $\alpha\in \Co_i$, one has
\begin{equation*}
\Phi_{r_{\alpha}}=\begin{cases} \Co_i,&\text{if }\alpha=\epsi_i,\\
(S(\alpha)\cup\{\alpha,2\epsi_i\})\setminus\{\epsi_i-\epsi_j\},&\text{if
}\alpha=\epsi_i+\epsi_j\text{ and }\Phi=C_n,\\
S(\alpha)\cup\{\alpha\}&\text{ otherwise}.
\end{cases}
\end{equation*}
By the way, $\Phi_{r_{\alpha}}$ consists of
$|S(\alpha)|+1=2|S^+(\alpha)|+1$ roots, so $s(\sigma)=1$,
$l(\sigma)=2|S^+(\beta)|+1$ and $\vartheta=0$. Thus,
$l(\sigma)-s(\sigma)-2\vartheta=2|S^+(\beta)|=2l$.

If $l>s_1$, then pick $i$ such that $s_1+\ldots+s_i<l\lee
s_1+\ldots+s_{i+1}$. As mentioned above, there exists $\beta\in
\Co_{2i+1}$ such that $|S^+(\beta)|=l-(s_1+\ldots+s_i)$. Set
$D=\{\beta_1,\ldots,\beta_i,\beta\}$. Then
$\Phi_{\sigma}=\cup_{j=1}^i\Phi_{r_{\beta_j}}\cup\Phi_{r_{\beta}}$
and these sets are disjoint. Hence,
\begin{equation*}
\begin{split}
l(\sigma)&=\sum_{j=1}^i|\Phi_{r_{\beta_j}}|+|\Phi_{r_{\beta}}|=
\sum_{j=1}^i(2|S^+(\beta_j)|+1)+(2|S^+(\beta)|+1)=\\
&=2(s_1+\ldots+s_i+|S^+(\beta)|)+(i+1)=2l+|D|
\end{split}
\end{equation*}
and $l(\sigma)-s(\sigma)-2\vartheta=2l+|D|-|D|-0=2l$. This concludes
the proof. $\square$

It follows from these results that if the ground field is
algebraically closed, then the dimension of~a~coadjoint orbit of the
group $U$ is equal to one of the numbers $0, 2,\ldots,2\mu$.

\bigskip In the remainder of the section we prove technical Lemmas
used in the proofs of Propositions~\ref{prop_pol_subal},
\ref{prop_polar} and \ref{prop_Fo_Fo_shtrih}. These Lemmas deal with
the case when $D$ contains \emph{two} roots from some column of
$\Phi^+$. Of~course, the proofs of these Lemmas are independent from
our previous results. \lemmp{Let $k$ be a field, $D\subset\Phi^+$ an
orthogonal subset, $\xi=(\xi_{\beta})_{\beta\in D}$ a set of
non-zero scalars from $k$. Let $\pt$, $\pt_0$, $\pt_2$, $\wt\pt$,
$\bt$ be defined as in the proof of Proposition
$\ref{prop_pol_subal}$. Then $[\pt_2+\wt\pt,\bt]\subset\pt$.
\label{lemm_p_0_1}}{Indeed, if $|D\cap \Co_1|\lee1$, then $\bt=0$.
Suppose $D\cap \Co_1=\{\epsi_1-\epsi_j,\epsi_1+\epsi_j\}$ (and,
consequently, $\Phi\neq C_n$ by Proposition \ref{prop_min_D} b)). In
this case, $\pt_1=\ut_1$ and $\at=0$.

Suppose that $x=\xi_{\epsi_1+\epsi_j}\cdot e_{\epsi_i-\epsi_j}-
\xi_{\epsi_1-\epsi_j}\cdot e_{\epsi_i+\epsi_j}\in\bt$, $1<i<j$, and
the coefficient of $e_{\gamma_0}$ in $y=\sum
y_{\gamma}e_{\gamma}\in\pt_2+\wt\pt$ is non-zero. Clearly,
$[x,e_{\gamma_0}]\neq0$ implies $\gamma_0\in \Ro_i$, because in our
case, $\pt\cap\sum_{\alpha\in \Co_j}ke_{\alpha}=0$ and $[x, x']=0$
for all $x'=\xi_{\epsi_1+\epsi_j}\cdot e_{\epsi_l-\epsi_j}-
\xi_{\epsi_1-\epsi_j}\cdot e_{\epsi_l+\epsi_j}$, $1<l<j$.

Let $\gamma_0=\epsi_l-\epsi_i$ for some $1<l<i$ (if $l=1$, then
$[x,e_{\gamma_0}]\in\pt_1$). It's easy to see that
$[x,e_{\gamma_0}]=cx'$, $c\in k^*$, where
$x'=\xi_{\epsi_1+\epsi_j}\cdot e_{\epsi_l-\epsi_j}-
\xi_{\epsi_1-\epsi_j}\cdot e_{\epsi_l+\epsi_j}$. But if
$x'\notin\bt$, then $D\cap \Ro_{-l}\neq\varnothing$, i.e.,
$\epsi_s+\epsi_l\in D$ for~some $1<s<l$ (see the definition of
$\pt_0$). If $\gamma_0\in\Mo_s$, then, by definition of $\pt_0$, the
coefficient of $e_{\gamma_0}$ in~$y$ is zero. Hence, $\gamma_0$ does
\emph{not} belong to $\Mo_s$.

But this means (see~(\ref{formula_minus})) that $\gamma_0$ or
$\epsi_s+\epsi_i=(\epsi_s+\epsi_l)-\gamma_0$ belongs to $\Mo_r$ for
some $1<r<s$, i.e., $D$ contains one of the roots of the form
$\epsi_r-\epsi_i$, $\epsi_r+\epsi_l$, $\epsi_r+\epsi_s$,
$\epsi_r+\epsi_i$. If $D$ contains one of the roots
$\epsi_r-\epsi_i$, $\epsi_r+\epsi_l$, then, by definition of
$\pt_0$, the coefficient of $e_{\gamma_0}$ in $y$ is zero. A root of
the form $\epsi_r+\epsi_s$ is not orthogonal to the root
$\epsi_s+\epsi_l\in D$, hence, $\epsi_r+\epsi_s$ does not belong to
$D$. Finally, if $\epsi_r+\epsi_i\in D$,\linebreak then
$D\cap\Ro_{-i}\neq0$, and, consequently, the vector $x$ doesn't
belong to $\pt_0$ (by definition of this subspace).

Thus, $[x,e_{\gamma_0}]\subset\pt$. Since $\gamma_0$ and $x$ are
arbitrary, $[y,\bt]\subset\pt$.}

\bigskip\lemmp{Let $k$ be a field, $D\subset\Phi^+$ an orthogonal
subset, $\xi=(\xi_{\beta})_{\beta\in D}$ a set of non-zero scalars
from $k$. Let $\pt$, $y$, $\gamma_0$, $\Co_i$ be defined as in the
proof of Proposition $\ref{prop_polar}$. Moreover, let $|D\cap
\Co_i|=2$. Then there exists $x'\in\pt+ky$ such that $f([y,
x'])\neq0$.\label{lemm_p_0_2}}{Let $D\cap \Co_i=\{\beta,\beta'\}$,
where $\beta=\epsi_i\pm\epsi_j$, $\beta'=\epsi_i\mp\epsi_j$, and
$\gamma_0\in \Co_l$, $i<l<j$. Put $x=\xi_{\epsi_i+\epsi_j}\cdot
e_{\epsi_l-\epsi_j}- \xi_{\epsi_i-\epsi_j}\cdot
e_{\epsi_l+\epsi_j}$. We assume without loss of generality that
$\beta=\epsi_i-\epsi_j$, $\beta'=\epsi_i+\epsi_j$ and
$\gamma_0=\epsi_l-\epsi_j$, $\gamma_0'=\epsi_l+\epsi_j$.

Since $\alpha_0=\beta-\gamma_0\in\Po$, Lemma \ref{lemm_one_sing} a)
shows that if $x$ and
$y_0=y_{\gamma_0}e_{\gamma_0}+y_{\gamma_0'}e_{\gamma_0'}$ are linear
independent, then $f([y,
e_{\alpha_0}])=\xi_{\beta}y_{\gamma_0}+\xi_{\beta'}y_{\gamma_0'}\neq0$,
so we can put $x'=e_{\alpha_0}$. On the other hand, if $y_0=cx$,
$c\in k$, and $x\in\pt_0$, then the coefficients of $e_{\gamma_0},
e_{\gamma_0'}$ in $y-cx\in\pt+ky$ are zero, so we can use induction
on the number of non-zero coefficients in $y$. Thus, it remains to
consider the case when $x$ and $y_0$ are linear dependent and
$x\notin\pt_0$.

This means that $D\cap \Ro_l\neq\varnothing$; in other words, there
exists $s<l$ such that $\epsi_s+\epsi_l\in D$. We claim that $s>i$.
Indeed, the roots $\gamma_0, \gamma_0'$ belong to $\Mo_i$, not to
$\Mo_s$, so if $s<i$, then there exists $r<s$ such that the roots
$\epsi_s\pm\epsi_j=(\epsi_s+\epsi_l)-(\epsi_l\mp\epsi_j)$ belong to
$\Mo_r$. But this contradicts the orthogonality~of~$D$ (see the
remark before Lemma \ref{lemm_one_sing}).

Hence, $i<s<l<j$. Consider the roots $\gamma_1=\epsi_s-\epsi_j$,
$\gamma_1'=\epsi_s+\epsi_j$. If one of them belongs to $\Mo_r$ for
some $r<i$, then the subset $D$ is not orthogonal, as in the case
when
$\epsi_i-\epsi_s={(\epsi_i\pm\epsi_j)}-{(\epsi_s\pm\epsi_j)}\in\Mo_r$
for some $r<i$. Hence, $\gamma_1,\gamma_1'\in\Mo_i$ (see the
definition of~$\Mo$). The vector $x'=\xi_{\epsi_i+\epsi_j}\cdot
e_{\gamma_1}- {\xi_{\epsi_i-\epsi_j}\cdot e_{\gamma_1'}}$ belongs to
$\pt_0$ (if $x'\notin\pt_0$, then $D\cap \Ro_s\neq\varnothing$,
which contradicts the orthogonality of $D$). Therefore,
$f([x,x'])=f(2\xi_{\beta}\xi_{\beta'}e_{\epsi_s+\epsi_l})=
2\cdot\xi_{\beta}\cdot\xi_{\beta'}\cdot\xi_{\epsi_s+\epsi_l}\neq0$.
Arguing as in Lemma~\ref{lemm_one_sing}, one can show that $f([y,
x'])=f([x, x'])\cdot y_{\gamma_0}/\xi_{\beta'}\neq0$.}

\bigskip\lemmp{Let $\Phi=B_n$ or $D_n$, and $D\cap
\Co_1=\{\epsi_1-\epsi_j,\epsi_1+\epsi_j\}$. Then
\begin{equation}
\begin{split}
l(\sigma)&=l'(\sigma')+|\Co_1|+|\Co_j|+4\cdot\#\{(l,
s)\mid1<l<s<j\text{ $\mathrm{and}$
}\epsi_l+\epsi_s\in\wt D\}+\\
&+2\cdot\#\{(l,s)\mid1<l<j<s\text{ $\mathrm{and}$
}\epsi_l-\epsi_s,\epsi_l+\epsi_s\in\wt
D\}+\\&+2\cdot\#\{l\mid1<l<j\text{ $\mathrm{and}$ }\epsi_l\in\wt
D\}.
\end{split}\label{formula_i_pm_j}
\end{equation}
\label{lemm_i_pm_j}}{As in Lemmas \ref{lemm_i_m_j} and
\ref{lemm_i_p_j}, $\Phi^+\setminus\wt\Phi^+=\Co_1\cup \Co_j\cup
\Ro_j\cup \Ro_{-j}$ (see.~(\ref{forumula_wt_Phi})). Clearly,
$\Co_1\subset\Phi_{\sigma}$. For~simplicity, put
$\sigma_1=r_{\epsi_1-\epsi_j}r_{\epsi_1+\epsi_j}$ (and so
$\sigma=\sigma_1\wt\sigma$). Let's study the action of
$\sigma_1$~and~$\wt\sigma$ on the roots from $\Co_j\cup \Ro_j\cup
\Ro_{-j}$. If $\alpha=\epsi_i\pm\epsi_j\in \Ro_j\cup \Ro_{-j}$, then
$\sigma_1(\alpha)=\epsi_i\mp\epsi_j>0$. If
$\alpha=\epsi_j\pm\epsi_l\in \Co_j$, then
$\sigma_1(\alpha)=-\epsi_j\pm\epsi_l<0$ (similarly,
$\sigma_1(\epsi_j)=-\epsi_j<0$). Thus,
$\Phi_{\sigma_1}\setminus\wt\Phi^+=\Co_1\cup \Co_j$. For the case
$\wt D=\varnothing$, there is nothing to prove.

Suppose $\beta=\epsi_l+\epsi_s\in\wt D$, where $1<l<s<j$. Clearly,
$r_{\beta}\sigma_1(\alpha)=\sigma_1(\alpha)<0$ for all $\alpha\in
\Co_j$. On~the~other hand, $r_{\beta}\sigma_1$ maps
$\epsi_l\pm\epsi_j$ and $\epsi_l\pm\epsi_j$ to the negative roots
$-\epsi_s\pm\epsi_j$ and $-\epsi_l\pm\epsi_j$ respectively. Hence,
$\Ro_j\cup \Ro_{-j}$ contains four roots with negative images under
the action of $r_{\beta}\sigma_1$. This~gives the~fourth summand in
the RHS of (\ref{formula_i_pm_j}).

Suppose $\beta=\epsi_l-\epsi_s,\beta'=\epsi_l+\epsi_s\in\wt D$,
where $1<l<j<s$. In this case, $r_{\beta}r_{\beta'}\sigma_1$ maps
$\epsi_l\pm\epsi_j\in \Ro_j\cup \Ro_{-j}$ and $\epsi_j\pm\epsi_s$ to
the negative roots $-\epsi_l\pm\epsi_j$ and $-\epsi_j\mp\epsi_s$
respectively. Note that the roots $\sigma_1(\epsi_j\pm\epsi_s)$ are
also negative, and
$r_{\beta}r_{\beta'}\sigma_1(\alpha)=\sigma_1(\alpha)$ for all other
$\alpha\in \Co_j\cup \Ro_j\cup \Ro_{-j}$. This gives the~fifth
summand in the RHS of (\ref{formula_i_pm_j}).

Now, suppose $\beta=\epsi_l\in\wt D$, where $1<j<l$, Then
$r_{\beta}\sigma_1$ maps $\epsi_l\pm\epsi_j\in \Ro_j\cup \Ro_{-j}$
to the negative roots $-\epsi_l\mp\epsi_j$, and
$r_{\beta}\sigma_1(\alpha)=\sigma_1(\alpha)$ for~all other
$\alpha\in \Co_j\cup \Ro_j\cup \Ro_{-j}$. This gives the last
summand in~the~RHS of (\ref{formula_i_pm_j}). Finally, suppose
$\beta=\epsi_l-\epsi_s\in\wt D$, where $1<l<j<s$. Then
$r_{\beta}\sigma_1$ maps $\epsi_l+\epsi_j\in \Ro_j$ to the negative
root $-\epsi_j+\epsi_s$, $r_{\beta}\sigma_1$ maps $\epsi_j+\epsi_s$
to the positive root $\epsi_l-\epsi_j$, and
$r_{\beta}\sigma_1(\alpha)=\sigma_1(\alpha)$ for~all other
$\alpha\in\Co_j\cup \Ro_j\cup \Ro_{-j}$. Thus, the number
$|\Phi_{\sigma}\setminus\wt\Phi^+|$ doesn't depend on roots from
$\wt D$ of~the~form~$\epsi_l-\epsi_s$. It's easy to see that the
action of $\sigma$ on $\Co_j\cup \Ro_j\cup \Ro_{-j}$ doesn't depend
on other roots from $\wt D$. This concludes the proof. }

\bigskip
Note that Lemmas \ref{lemm_p_0_1}, \ref{lemm_p_0_2}, Propositions
\ref{prop_pol_subal}, \ref{prop_polar} and Theorem \ref{mtheo_pol}
are also valid for a field $k$ of \emph{zero} characteristic
(indeed, their proofs do \emph{not} depend on the characteristic of
the ground field). In~particular this allows to find polarizations
for orbit associated with orthogonal subsets for~the~case $k=\Rp$
(they play an important role in the construction of unitary
irreducible representations of~corresponding nilpotent Lie groups,
see, for example, \cite[p. 182]{Kirillov1}).

\end{document}